\documentclass[a4paper,11pt,reqno]{amsart}

\usepackage{amscd}

 \font\sevenrm=cmr7

\input prepictex
\input pictex
\input postpictex

\newtheorem{theorem}[equation]{Theorem}
\newtheorem{lemma}[equation]{Lemma}
\newtheorem{corollary}[equation]{Corollary}
\newtheorem{proposition}[equation]{Proposition}
\numberwithin{equation}{section}

\theoremstyle{definition}
\newtheorem*{example}{Example}
\newtheorem{remark}[equation]{Remark}
\newtheorem*{acknow}{Acknowledgments}

\newcommand{\N}{{\mathbb N}}
\newcommand{\Z}{{\mathbb Z}}
\newcommand{\R}{{\mathbb R}}
\newcommand{\C}{{\mathbb C}}
\newcommand{\J}{{\mathcal J}}
\newcommand{\cP}{{\mathcal P}}
\newcommand{\g}{{\mathfrak g}}
\newcommand{\h}{{\mathfrak h}}

 \DeclareMathOperator{\Br}{Br}
 \DeclareMathOperator{\End}{End}
 \DeclareMathOperator{\Hom}{Hom}
 \DeclareMathOperator{\alg}{alg}
 \DeclareMathOperator{\diag}{diag}
 \DeclareMathOperator{\Real}{Re}
 \DeclareMathOperator{\Dr}{D_r^{\text{\rm marked}}}

\newcommand{\OrV}{\otimes_\R^rV}
\newcommand{\OrVd}{\otimes_\R^rV_\R^*}
\newcommand{\OdrVd}{\otimes_\R^{2r}V_\R^*}
\newcommand{\EUVr}{\End_{U(V,h)}(\OrV)}
\newcommand{\uvh}{\mathfrak{u}(V,h)}
\newcommand{\suvh}{\mathfrak{su}(V,h)}
\newcommand{\glv}{\mathfrak{gl}_\C(V)}
\newcommand{\slv}{\mathfrak{sl}_\C(V)}
\newcommand{\gln}{\mathfrak{gl}_n}

\def\inicioYoung{
       \begingroup
       \def\vr{\vrule height 10pt width .5pt depth 2pt}
       \def\caja##1{\lower 2pt\vbox{\offinterlineskip
                    \hrule height .5pt
                    \hbox to 14pt{\vr\hfill##1\hfill\vr}
                    \hrule height .5pt}}
       \vtop\bgroup \offinterlineskip \tabskip=-.5pt \lineskip=-.5pt
            \halign\bgroup &\caja{##\unskip}\unskip  \crcr }
\def\finYoung{\egroup\egroup\hskip 3pt\endgroup}

\newenvironment{Young}{\inicioYoung}{\finYoung}

\begin{document}

\title[Modified Brauer algebra]{A modified Brauer algebra as centralizer algebra of the
unitary group}

\author{Alberto Elduque}

\thanks{Supported by the Spanish Ministerio de Ciencia y
Tecnolog\'{\i}a and FEDER (BFM 2001-3239-C03-03)}

\address{Departamento de Matem\'aticas, Universidad de
Zaragoza, 50009 Zaragoza, Spain}

\email{elduque@unizar.es}

\date{\today}



\begin{abstract}
The centralizer algebra of the action of $U(n)$ on the real tensor
powers $\OrV$ of its natural module, $V=\C^n$, is described by
means of a modification in the multiplication of the signed Brauer
algebras. The relationships of this algebra with the invariants
for $U(n)$ and with the decomposition of $\OrV$ into irreducible
submodules is considered.
\end{abstract}

\maketitle


\section{Introduction}

The motivation for this work comes from a paper by Gray and
Hervella \cite{Gray-Hervella}: Let $(M,g,J)$ be an almost
Hermitian manifold; that is, $M$ is a Riemannian manifold with
Riemannian metrig $g$, and endowed with an almost complex
structure $J$. Let $\nabla$ be the Riemannian connection and $F$
the K\"ahler form: $F(X,Y)=g(JX,Y)$ for any $X,Y\in \chi(M)$ (the
set of smooth vector fields). The tensor $G=\nabla F$ satisfies
$G(X,Y,Z)=-G(X,Z,Y)=-G(X,JY,JZ)$ for any $X,Y,Z\in \chi(M)$.
Therefore, at any point $p\in M$, $\alpha=G_p$ belongs to
\begin{multline}\label{e:Wp}
W_p=\{ \alpha\in M_p^*\otimes_{\R}M_p^*\otimes_\R M_p^* :\\
 \alpha(x,y,z)=-\alpha(x,z,y)=-\alpha(x,Jy,Jz)\ \forall x,y,z\in
 M_p\},
\end{multline}
where $M_p$ denotes the tangent space at $p$ and $M_p^*$ its dual
(the cotangent space), and $M_p^*\otimes_{\R}M_p^*\otimes_\R
M_p^*$ is identified naturally with the space of trilinear forms
on $M_p$.

The classification of almost hermitian manifolds in
\cite{Gray-Hervella} is based on the decomposition of $W_p$ into
irreducible modules under the action of the unitary group $U(n)$,
and this is done by first providing four specific subspaces of
$W_p$ (if the dimension of $M$ is not very small) and then showing
that they are irreducible (by means of invariants) and $W_p$ is
their direct sum. No clue is given about how these four subspaces
are obtained. A different way to obtain this decomposition is
given in \cite{Pedroetal} based on complexification of $W_p$ and
the use of Young symmetrizers.

The situation above extends naturally to the following problem:

\medskip

\emph{Given a complex vector space $V$ of dimension $n$, endowed
with a nondegenerate hermitian form $h:V\times V\rightarrow \C$,
decompose the $n^{\text{th}}$ tensor power  $\OrV$ (over the real
numbers!) into a direct sum of irreducible modules for the unitary
group $U(V,h)=\{g\in GL_\C(V) : h(gv,gw)=h(v,w)\ \forall v,w\in
V\}$.}

\medskip

Here, the convention is that $h(\alpha v,w)=\alpha h(v,w)$ and
$h(v,w)=\overline{h(w,v)}$ for any $\alpha\in \C$, $v,w\in V$.

\medskip

The well-known Schur-Weyl duality \cite{Schur1,Schur2,Weyl}
relates the representation theory of the general linear group
$GL_\C(V)$ with that of the symmetric group $S_r$ via the
naturally centralizing actions of the two groups on the space
$\otimes_\C^rV$: $GL_\C(V)\rightarrow \otimes_\C^rV \leftarrow
S_r$. Brauer \cite{Brauer} considered the analogous situation for
the orthogonal and symplectic groups, where $S_r$ has to be
replaced by what are now called the Brauer algebras:
$O(V)\rightarrow \otimes_\C^rV\leftarrow \Br_r(n)$ and
$Sp(V)\rightarrow \otimes_\C^rV\leftarrow \Br_r(-n)$. More
recently, Brauer algebras and their generalizations, specially the
BMW algebra, have been looked at in the context of quantum groups
and low dimensional topology \cite{Jimbo,BW,Mura,HRam,LRam}.

\medskip

In our problem, the decomposition of $\otimes_\R^rV$ into a direct
sum of irreducible modules for $U(V,h)$ is intimately related to
the action of the centralizer algebra
$\End_{U(V,h)}(\otimes_\R^rV)$, and the main part of the paper
will be devoted to computing this centralizer algebra. This will
be done, following a classical approach, by relating it to the
multilinear $U(V,h)$-invariant maps
$f:V\times\stackrel{r}{\cdots}\times V\rightarrow \R$. These
invariants will be the subject of Section 2. Section 3 will be
devoted to the determination of the centralizer algebra, while
Section 4 will give a combinatorial description of it, as well as
a presentation by generators and relations. It will turn out that
the centralizer algebra looks like the Signed Brauer Algebra
considered in \cite{Par1,Par2} (see Remark \ref{r:signed}). This
algebra appears, for sufficiently large dimension, as the
centralizer algebra of the action of the product of orthogonal
groups $O\left(S^2(V)\right)\times O\left(\Lambda^2(V)\right)$ on
$\otimes_\C^r(V\otimes_\C V)=\otimes_\C^r\left(S^2(V)\oplus
\Lambda^2(V)\right)$, for a vector space $V$ equipped with a
nondegenerate symmetric bilinear form $b:V\times V\rightarrow \C$,
which induces a nondegenerate bilinear form on $V\otimes_\C
V=S^2(V)\oplus \Lambda^2(V)$ (orthogonal direct sum). In section 5
it will be shown how to use the information on the centralizer
algebra, together with the results in \cite{Benkartetal}, to
decompose $\otimes_\R^rV$ into a direct sum of irreducible
$U(V,h)$-modules. A couple of examples will be given: the one in
\cite{Gray-Hervella} mentioned above, and another one considered
in \cite{Abbena-Garbiero}, used to classify homogeneous K\"ahler
structures.

\section{Invariants}

This section is devoted to prove the next result:

\begin{theorem}\label{t:invariantes}
Let $V$ be an $n$-dimensional complex vector space, endowed with a
nondegenerate hermitian form $h:V\times V\rightarrow \C$ and let
$r\in \N$. If $f:V\times\stackrel{r}{\cdots}\times V\rightarrow
\R$ is a nonzero multilinear $U(V,h)$-invariant form, then $r$ is
even ($r=2m$) and $f$ is a linear combination of the invariant
maps:
\begin{equation*}
\begin{split}
V\times\stackrel{2m}{\cdots}\times V&\rightarrow \R\\
(v_1,\ldots,v_{2m})&\mapsto \prod_{l=1}^m \left\langle
v_{\sigma(2l-1)},J^{\delta_l}v_{\sigma(2l)}\right\rangle
\end{split}
\end{equation*}
where $\sigma\in S_{2m}$ (the symmetric group on
$\{1,\ldots,2m\}$), $\delta_1,\ldots,\delta_m\in \{0,1\}$,
$J:V\rightarrow V$ is the multiplication by $i\in \C$ ($i^2=-1$),
and $\langle\,\mid\,\rangle$ denotes the real part of $h$ (so that
$h(v,w)=\langle v\mid w\rangle +i\langle v\mid Jw\rangle$ for any
$v,w\in V$).
\end{theorem}

In case $\dim_\C V\geq r$, this appears in \cite{Iwahori}. For
arbitrary $r$, it is asserted in \cite{Gray-Hervella} without
proof. A proof will be provided here, which will be based on
methods to be used later on.

\bigbreak

\emph{Throughout the paper $(V,h)$, $J$ and
$\langle\,\mid\,\rangle$ will be assumed to satisfy the hypotheses
of Theorem \ref{t:invariantes}.}

\bigbreak

Let $r\in \N$, for any $l\in \{1,\ldots,r\}$ consider the
$\R$-linear map:
\begin{equation*}
\begin{split}
J_l:\ \otimes_\R^r V\ &\longrightarrow\ \otimes_\R^rV\\
v_1\otimes\cdots\otimes v_r&\mapsto v_1\otimes\cdots\otimes
Jv_l\otimes\cdots\otimes v_r
\end{split}
\end{equation*}
(action of $J$ on the $l^{\text{th}}$-spot). Then $J_l\in
\End_{U(V,h)}(\OrV)$, the centralizer algebra.

\smallskip

\emph{As a general rule, the elements of the centralizer algebra
will act on the right.}

\smallskip

Let $\J=\alg_\R\{ J_1,\ldots,J_r\}$ the (real) subalgebra of
$\EUVr$ generated by the $J_l$'s. It is clear the $\J$ is
isomorphic, as an algebra, to $\otimes_\R^r\C$, under the map that
sends $J_l$ to $1\otimes\cdots\otimes i\otimes\cdots\otimes 1$
($i$ in the $l^{\text{th}}$ slot) for any $l$. Note that for any
$1\leq l\ne m\leq r$, $\frac{1}{2}(1\pm J_lJ_m)$ is an idempotent
in $\J$ since $J_l^2=-1$ for any $l$. For any nonempty subset
$\cP\subseteq \{1,\ldots,r\}$ and any $p\in \cP$, let
$\cP^c=\{1,\ldots,r\}\setminus \cP$ and consider the following
element of $\J$:
$$
e_\cP =\frac{1}{2^{r-1}}\prod_{p\ne q\in \cP}\hskip -3pt
(1-J_pJ_q) \prod_{q\in \cP^c}(1+J_pJ_q)\,.
$$
Then:

\begin{proposition}\label{p:2.2}
Under the conditions above:
\begin{enumerate}
\item $e_\cP$ does not depend on the chosen element $p\in\cP$.
\item $e_\cP$ is a primitive idempotent of $\J$.
\item Given any $p\in \{1,\ldots,r\}$, $\J=\oplus_{p\in
\cP\subseteq\{1,\ldots,r\}}\C e_\cP$.
\end{enumerate}
\end{proposition}

\begin{proof}
The $\R$-linear map $\C\otimes_\R\C\rightarrow \C\oplus\C$,
$\alpha\otimes\beta\mapsto (\alpha\beta,\alpha\bar\beta)$, yields
an algebra isomorphism. Therefore, as real algebras, $\J\cong
\otimes_\R^r\C\cong \C^{2^{r-1}}$. Now, fix $p\in \cP\subseteq
\{1,\ldots,r\}$; then if $p\in\cP'\subseteq \{1,\ldots,r\}$ and
$q\in \cP\setminus\cP'$, $e_\cP e_{\cP'}$ contains the factor
$(1-J_pJ_q)(1+J_pJ_q)=0$, and hence $e_\cP e_{\cP'}=0$. The same
argument works for any $q\in \cP'\setminus \cP$. Therefore $e_\cP$
and $e_{\cP'}$ are orthogonal idempotents. Since there are
$2^{r-1}$ subsets $\cP\subseteq\{1,\ldots,r\}$ containing $p$ and
$\J$ is an algebra over $\C$, where the action of $\C$ is given
``on the $p^{\text{th}}$ slot'' ($\J\cong\otimes_\R^r\C$), to
prove (2) and (3) it is enough to check that $e_{\cP}$ is nonzero
for any $p\in\cP\subseteq \{1,\ldots,r\}$. For simplicity, and
without loss of generality, assume $p=1$. Then $2^{r-1}e_\cP$ is
the sum of $2^{r-1}$ summands $\pm J_1^mJ_2^{\delta_2}\cdots
J_r^{\delta_r}$, with $\delta_l=0$ or $1$, $m\geq 0$ and
$m=\delta_2+\cdots+\delta_r$. All these summands are linearly
independent (over $\R$) in $\J\cong\otimes_\R^r\C$.

It remains (1) to be proved. Take $p\ne p'\in \cP$, $e_\cP
=\frac{1}{2^{r-1}}\prod_{p\ne q\in\cP}(1-J_pJ_q)\prod_{q\in
\cP^c}(1+J_pJ_q)$ and $e'_\cP =\frac{1}{2^{r-1}}\prod_{p'\ne
q\in\cP}(1-J_{p'}J_q)\prod_{q\in \cP^c}(1+J_{p'}J_q)$. Notice that
for any $s\ne t$ in $\{1,\ldots,r\}$
\begin{gather*}
(1-J_sJ_t)J_s=J_s+J_t=(1-J_sJ_t)J_t\\
(1+J_sJ_t)J_s=J_s-J_t=-(1+J_sJ_t)J_t
\end{gather*}
and hence
\begin{equation}\label{eqProp2.2}
\begin{cases}
e_\cP J_p=e_\cP J_q &\text{for any $q\in\cP$,}\\
e_\cP J_p=-e_\cP J_q &\text{for any $q\in \cP^c$.}
\end{cases}
\end{equation}
So that, since $J_p^2=-1$,
\begin{equation*}
\begin{split}
e_\cP(1-J_sJ_t)&=\begin{cases} 0&\text{if either $s\in\cP$,
$t\in\cP^c$ or $s\in\cP^c$, $t\in\cP$,}\\
2e_\cP &\text{if either $s,t\in\cP$ or $s,t\in\cP^c$,}
\end{cases}
\intertext{while}
e_\cP(1+J_sJ_t)&=\begin{cases}
 0&\text{if either $s,t\in\cP$ or $s,t\in\cP^c$,}\\
2e_\cP &\text{if either $s\in\cP$, $t\in\cP^c$ or $s\in\cP^c$,
$t\in\cP$,}
\end{cases}
\end{split}
\end{equation*}
Therefore $e_\cP e'_\cP=e_\cP$ and, with the same argument, $e_\cP
e'_\cP=e'_\cP$, whence $e_\cP=e'_\cP$.
\end{proof}

\begin{corollary}\label{c:2.3}
For any fixed $1\leq p\leq r$, $1=\sum_{p\in\cP\in\{1,\ldots,r\}}
e_\cP$. \qed
\end{corollary}

Also, equation \eqref{eqProp2.2} immediately yields:

\begin{lemma}\label{l:2.4}
Let $\emptyset\ne\cP\subseteq\{1,\ldots,r\}$ and $p\in \cP$, then
\begin{equation*}
\quad\bigl(\OrV\bigr)e_\cP=\{x\in\OrV: xJ_q=xJ_p\ \forall
q\in\cP,\
 xJ_q=-xJ_p\ \forall q\in\cP^c\} .\quad \qed
\end{equation*}
\end{lemma}

Now, denote by $V^*$ the dual vector space of $V$ as a complex
vector space. $V^*$ is a module too for the unitary group
$U(V,h)$. Take $\emptyset \ne \cP\subseteq\{1,\ldots,r\}$ and
consider:
\begin{equation*}
V_l=\begin{cases}
 V&\text{if $l\in \cP$,}\\
 V^*&\text{if $l\in \cP^c$.}
 \end{cases}
\end{equation*}

Then:

\begin{proposition}\label{p:2.5}
The $\R$-linear map
\begin{equation*}
\begin{split} \Phi_\cP :\ \bigl(\OrV\bigr)e_\cP\ &\longrightarrow
  V_1\otimes_\C V_2\otimes_\C \cdots\otimes_\C V_r\\
  (v_1\otimes\cdots\otimes v_r)e_\cP&\mapsto \quad
  w_1\otimes w_2\otimes\cdots\otimes w_r
\end{split}
\end{equation*}
where $w_l=v_l$ if $l\in\cP$ and $w_l=h(-,v_l)\in V^*$ if
$l\in\cP^c$, is well defined and an isomorphism of
$U(V,h)$-modules.
\end{proposition}
\begin{proof}
The linear map
\begin{equation}\label{eqProp2.5}
\begin{split} \Psi_\cP :\ \OrV\ &\longrightarrow
  V_1\otimes_\C V_2\otimes_\C \cdots\otimes_\C V_r\\
  v_1\otimes\cdots\otimes v_r&\mapsto \quad
  w_1\otimes w_2\otimes\cdots\otimes w_r
\end{split}
\end{equation}
with $w_1,\ldots,w_r$ as above, is well defined and a homomorphism
of $U(V,h)$-modules. Besides, if $p,q\in\cP$, then $w_p=v_p$ and
$w_q=v_q$ above, so
\begin{equation*}
\begin{split}
\Psi_\cP\Bigl((v_1&\otimes\cdots\otimes
v_r)\frac{1}{2}(1-J_pJ_q)\Bigr)\\
&=\frac{1}{2}\Psi_\cP\bigl(v_1\otimes\cdots\otimes v_r-
      v_1\otimes\cdots\otimes iv_p\otimes\cdots\otimes iv_q\otimes
       \cdots\otimes v_r\bigr)\\
&=\frac{1}{2}\bigl(w_1\otimes\cdots\otimes w_r-
      w_1\otimes\cdots\otimes iw_p\otimes\cdots\otimes iw_q\otimes
       \cdots\otimes w_r\bigr)\\
&=w_1\otimes\cdots\otimes w_r,
\end{split}
\end{equation*}
while if $p\in\cP$ and $q\in\cP^c$, then $w_p=v_p$ and
$w_q=h(-,v_q)$ so, since $h(-,iv_q)=-ih(-,v_q)=-iw_q$,
\begin{equation*}
\begin{split}
\Psi_\cP\Bigl((v_1&\otimes\cdots\otimes
v_r)\frac{1}{2}(1+J_pJ_q)\Bigr)\\
&=\frac{1}{2}\Psi_\cP \bigl(v_1\otimes\cdots\otimes v_r+
      v_1\otimes\cdots\otimes iv_p\otimes\cdots\otimes iv_q\otimes
       \cdots\otimes v_r\bigr)\\
&=\frac{1}{2}\bigl(w_1\otimes\cdots\otimes w_r-
      w_1\otimes\cdots\otimes iw_p\otimes\cdots\otimes iw_q\otimes
       \cdots\otimes w_r\bigr)\\
&=w_1\otimes\cdots\otimes w_r.
\end{split}
\end{equation*}
Therefore, since $e_\cP=\prod_{p\ne
q\in\cP}\left(\frac{1}{2}(1-J_pJ_q)\right)
\prod_{q\in\cP^c}\left(\frac{1}{2}(1+J_pJ_q)\right)$, it follows
that $\Psi_\cP (x)=\Psi_\cP (xe_\cP)$ for any $x\in \OrV$, and
then $\Psi_\cP$ restricts to $\Phi_\cP$, which is thus well
defined. The inverse is given by
\begin{equation*}
\begin{split}
\Phi_\cP^{-1}: V_1\otimes_\C \cdots\otimes_\C V_r&\longrightarrow
\bigl(\OrV\bigr) e_\cP\\
w_1\otimes\cdots\otimes w_r\ &\mapsto (v_1\otimes \cdots\otimes
v_p)e_\cP
\end{split}
\end{equation*}
where $v_l=w_l$ if $l\in \cP$, while $w_l=h(-,v_l)$ for a unique
$v_l\in V$ if $l\in\cP^c$. This is well defined because of Lemma
\ref{l:2.4}
\end{proof}

Now, Corollary \ref{c:2.3} and Proposition \ref{p:2.5} yield:

\begin{corollary}\label{c:2.6}
Fix $p\in\{1,\ldots,r\}$,  then the $U(V,h)$-module $\OrV$ is
isomorphic to
$$
\bigoplus_{p\in\cP\subseteq\{1,\ldots,r\}}\hskip -10pt
V_{1\cP}\otimes_\C\cdots\otimes_\C V_{r\cP}
$$
where $V_{l\cP}=V$ if $l\in\cP$, while $V_{l\cP}=V^*$ otherwise.
\qed
\end{corollary}

Notice that $\OrV$ is a complex vector space with the action of
$\C$ on the $p^{\text{th}}$ slot, and that the isomorphism in
Corollary \ref{c:2.6} is then an isomorphism of complex vector
spaces too.

\smallskip

The final prerequisite in the proof of Theorem \ref{t:invariantes}
is the next straightforward result:

\begin{lemma}\label{l:2.7}
Let $\g$ be a real Lie algebra, $\rho:\g\rightarrow \End_\C(W)$ a
complex representation of $\g$, $f:W\rightarrow \R$ a linear
$\g$-invariant map ($f(x.w)=f(w)$ for any $x\in\g$, $w\in W$).
Then there is a complex linear $\g$-invariant map $g:W\rightarrow
\C$ such that $f$ is the real part of $g$ ($f=\Real g$). \qed
\end{lemma}

\bigbreak

\begin{proof}[Proof of Theorem \ref{t:invariantes}]
 From the previous results we obtain:
{\small
\begin{equation*}
\begin{split}
\{ &S:  V\times\stackrel{r}{\cdots}\times V\rightarrow \R \mid S\
\text{multilinear and $U(V,h)$-invariant}\}\\
 &\simeq\{ S:\OrV\rightarrow \R \mid S\ \text{linear and
  $U(V,h)$-invariant}\}\\
 &=\hskip -15pt\bigoplus_{1\in\cP\subseteq\{1,\ldots,r\}}\hskip -10pt\{ S_\cP
  :\bigl(\OrV\bigr)e_\cP\rightarrow \R \mid S_\cP\ \text{linear
  and $U(V,h)$-invariant}\}\\
 &\simeq\hskip -15pt\bigoplus_{1\in\cP\subseteq\{1,\ldots,r\}}\hskip -10pt
  \{ S_\cP
  :V_{1\cP}\otimes_\C\cdots\otimes_\C V_{r\cP}\rightarrow \R \mid S_\cP\
  \text{linear and $U(V,h)$-invariant}\}\\
 &=\hskip -15pt\bigoplus_{1\in\cP\subseteq\{1,\ldots,r\}}\hskip -10pt
  \{ \Real{T_\cP} \mid
   T_{\cP}:V_{1\cP}\otimes_\C\cdots\otimes_\C V_{r\cP}\rightarrow\C
  \text{ ($\C$)-linear and $U(V,h)$-invariant}\}\\
 &=\hskip -15pt\bigoplus_{1\in\cP\subseteq\{1,\ldots,r\}}\hskip -10pt
   \{ \Real{T_\cP} \mid
   T_{\cP}:V_{1\cP}\otimes_\C\cdots\otimes_\C V_{r\cP}\rightarrow\C
  \text{ ($\C$)-linear and $GL_\C(V)$-invariant}\}
\end{split}
\end{equation*}
}(The last equality is due to the fact that $U(V,h)$ is a form of
$GL_\C(V)$.) But the invariant theory of $GL_\C(V)$ \cite{Weyl}
shows that
$$
\{T :V_{1\cP}\otimes_\C\cdots\otimes_\C V_{r\cP}\rightarrow\C \mid
  \text{$T$ is $\C$-linear and $GL_\C(V)$-invariant.}\}
$$
is trivial unless $r$ is even, $r=2m$, and $\cP$ contains exactly
$m$ elements. In this latter case, $\cP=\{l_1,\ldots,l_m\}$
($l_1=1$), $\cP^c=\{s_1,\ldots,s_m\}$ and any such invariant $T$
is a (complex) linear combination of invariants of the form
\begin{equation*}
w_1\otimes \cdots\otimes w_{2m}\mapsto \prod_{j=1}^m
\varphi_{s_{\sigma(j)}}(v_{l_j})
\end{equation*}
where $w_l=v_l\in V$ for $l\in\cP$, $w_l=\varphi_l\in V^*$ for
$l\in\cP^c$, and $\sigma\in S_m$, the symmetric group on
$\{1,\ldots,m\}$.

Taking into account the definitions of the isomorphisms
$\Phi_{\cP}$ and homomorphisms $\Psi_\cP$ in Proposition
\ref{p:2.5} and Equation \eqref{eqProp2.5}, any $U(V,h)$-invariant
linear map $T:\OrV\rightarrow \C$ is a complex linear combination
of the maps
\begin{equation*}
v_1\otimes\cdots\otimes v_{2m}\mapsto \prod_{j=1}^m
h\bigl(v_{\sigma(2j-1)},v_{\sigma(2j)}\bigr)
\end{equation*}
where $\sigma\in S_{2m}$. Since $h(v,w)=\langle v\mid w\rangle +
i\langle v\mid Jw\rangle$, Lemma \ref{l:2.7} finishes the proof.
\end{proof}

\bigbreak

A final remark for this section is that using the invariant theory
for $SL_\C(V)$ instead of $GL_\C(V)$ and the same arguments as
above, one arrives at:

\begin{proposition}
The invariant multilinear $SU(V,h)$-invariant maps
$f:V\times\stackrel{r}{\cdots}\times V\rightarrow \R$ are exactly
the linear combinations of the maps:
\begin{multline}
(v_1,\ldots,v_r)\mapsto\\
 \Bigl(\prod_{j=1}^m \left\langle
  v_{\sigma(2j-1)},v_{\sigma(2j)}\right\rangle\Bigr)
 \Bigl(\prod_{l=0}^{s-1}
 \det{}_\C\bigl(v_{\sigma(2m+nl+1},\cdots,
    v_{\sigma(2m+n(l+1))}\bigr)\Bigr)
\end{multline}
where $n=\dim_\C V$, $m,s\geq 0$ with $n=2m+ns$ and $\sigma\in
S_r$. \qed
\end{proposition}

\section{Centralizer algebra}

To compute the centralizer algebra $\End_{U(V,h)}(\OrV)$, it is
enough to use the fact that $\End_\R(\OrV)$ is isomorphic (as
vector spaces and as $U(V,h)$-modules) to
$\bigl(\OrVd\bigr)\otimes_\R\bigl(\OrV\bigr)$, where $V_\R^*$
denotes the dual as a real vector space, to distinguish it from
$V^*=\Hom_\C(V,\C)$ (the dual as a complex vector space). But $V$
is isomorphic to $V_\R^*$ as $U(V,h)$-module by means of
$\langle\,\mid\,\rangle$ ($V\rightarrow V_\R^*$, $v\mapsto \langle
v\mid -\rangle$, is an isomorphism). Therefore, $\End_\R(\OrV)$ is
isomorphic to $\OdrVd$, which is naturally identified with the
space of multilinear maps: $f: V\times\stackrel{2r}{\cdots}\times
V\rightarrow\R$. Under these isomorphims, the centralizer algebra
$\End_{U(V,h)}(\OrV)$ (which is the subalgebra of $\End_\R(\OrV)$
fixed under the action of $U(V,h)$) corresponds to the space of
multilinear and $U(V,h)$-invariant maps.

Hence, to compute $\End_{U(V,h)}(\OrV)$ one has just to keep track
of the isomorphisms above. Let us proceed with an example:
consider the multilinear $U(V,h)$-invariant map
\begin{equation}\label{e:ejemplo-f}
\begin{split}
f: V\times V\times V\times V\times V\times V&\longrightarrow\ \R\\
(v_1,v_2,v_3,v_4,v_5,v_6)\ &\mapsto
 \langle v_1\mid Jv_3\rangle\langle v_2\mid Jv_4\rangle
  \langle v_5\mid Jv_6\rangle
\end{split}
\end{equation}
and let $\{e_l\}_{l=1}^{2n}$ be a basis (over $\R$) of $V$, and
$\{f_l\}_{l=1}^{2n}$ its dual basis relative to
$\langle\,\mid\,\rangle$ (so that $\langle e_p\mid f_q\rangle
=\delta_{pq}$ for any $p,q$). Let $e_l^*=\langle e_l\mid
-\rangle$, $f_l^*=\langle f_l\mid -\rangle\in V_\R^*$ for any $l$.
Then notice that the bilinear $U(V,h)$-invariant maps
$(v,w)\mapsto \langle v\mid w\rangle$ and $(v,w)\mapsto \langle
v\mid Jw\rangle$ correspond in $V_\R^*\otimes_\R V_\R^*$ to
$\sum_{l=1}^{2n} e_l^*\otimes f_l^*$ and
$\sum_{l=1}^{2n}e_l^*\otimes (f_l^*\circ J)$ respectively. Thus,
the multilinear map $f$ in \eqref{e:ejemplo-f} corresponds to:
\begin{equation}\label{e:identificaciones}
\begin{split}
f&\simeq \sum_{a,b,c=1}^{2n} e_a^*\otimes e_b^*\otimes (f_a^*\circ
   J)\otimes (f_b^*\circ J)\otimes e_c^*\otimes (f_c^*\circ J)
\text{\qquad in $\otimes_\R^6 V_\R^*$}\\
 &\simeq (-1)^2\sum_{a,b,c=1}^{2n} e_a^*\otimes e_b^*\otimes (f_a^*\circ
   J)\otimes Jf_b\otimes e_c\otimes Jf_c\\
\intertext{\qquad\qquad\qquad\qquad in
  $(\otimes_\R^3V_\R^*)\otimes_\R(\otimes_\R^3V)$ (since $\langle
  f_a\mid Jv\rangle =-\langle Jf_a\mid v\rangle$)}
 &\simeq \left( v_1\otimes v_2\otimes v_3\mapsto
   \sum_{a,b,c=1}^{2n} e_a^*(v_1)e_b^*(v_2)f_a^*(Jv_3)\,
      Jf_b\otimes e_c\otimes Jf_c\right)\\
\intertext{\hskip 3.5in in $\End_\R(\otimes_\R^3 V)$}
 &= \left( v_1\otimes v_2\otimes v_3\mapsto
    \langle v_1\mid Jv_3\rangle \sum_{c=1}^{2n}
     \bigl(\sum_{b=1}^{2n}e_b^*(v_2)Jf_b\bigr)\otimes e_c\otimes
     Jf_c\right)\\
 &= \left( v_1\otimes v_2\otimes v_3\mapsto
    \langle v_1\mid Jv_3\rangle \sum_{c=1}^{2n}
     Jv_2\otimes e_c\otimes Jf_c\right)\\
 &=J_3c_{13}J_2J_3(12)\in\End_{U(V,h)}(\otimes_\R^3V)
\end{split}
\end{equation}
where
$$
(v_1\otimes v_2\otimes v_3)c_{13}=\langle v_1\mid v_3\rangle
\sum_{a=1}^{2n} e_a\otimes v_2\otimes f_a
$$
and $(12)$ denotes the permutation of the first two slots. Notice
that $c_{13}$ does not depend on the chosen bases.

As we have seen, the $J_l$'s belong to the centralizer algebra
$\End_{U(V,h)}(\OrV)$, and so do the \emph{contraction maps}
$c_{pq}$ ($1\leq p<q\leq r$) defined as above:
\begin{equation}\label{e:contracciones}
\begin{split}
&(v_1\otimes\cdots\otimes v_r)c_{pq}\\
&= \langle v_p\mid v_q\rangle
 \sum_{a=1}^{2n} v_1\otimes\cdots\otimes v_{p-1}\otimes e_a\otimes
 v_{p+1}\otimes\cdots\otimes v_{q-1}\otimes f_a\otimes
 v_{q+1}\cdots\otimes v_r
\end{split}
\end{equation}

The arguments used for this particular $f$ in \eqref{e:ejemplo-f}
work in general. Therefore:

\begin{theorem}\label{t:centralizador}
The centralizer algebra $\End_{U(V,h)}(\OrV)$ is generated (as a
real algebra) by the $J_l$'s, $c_{pq}$'s and the action of the
symmetric group $S_r$:
\begin{equation*}
\ \quad End_{U(V,h)}(\OrV)=\alg_\R\bigl\{\rho(S_r),J_l,c_{pq} :
l,p,q=1,\ldots,r,\, p<q\bigr\}.\qquad \qed
\end{equation*}
\end{theorem}

\smallskip

(Given $\sigma\in S_r$, $\rho(\sigma)$ denotes the map
$v_1\otimes\ldots\otimes v_r\mapsto v_{\sigma(1)}\otimes\cdots
\otimes v_{\sigma(r)}$.)

\bigskip

\section{Combinatorial description}

Consider the element in \eqref{e:identificaciones}, which belongs
to the centralizer algebra $\End_{U(V,h)}(\otimes_\R^6V)$:
$$
v_1\otimes v_2\otimes v_3\mapsto \langle v_1\mid Jv_3\rangle
 \sum_{a=1}^{2n}Jv_2\otimes e_a\otimes Jf_a.
$$
(Notation as in the previous section.) It will be represented by
the \emph{marked diagram}:
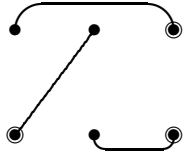
\begin{equation*}
\beginpicture
 \setcoordinatesystem
 \setplotarea x from 0 to 60, y from -40 to 40
 \linethickness=.8pt
 \setplotsymbol ({\sevenrm .})
 \multiput{$\bullet$} at 0 20 30 20 60 20 /
 \multiput{$\bullet$} at 0 -20 30 -20 60 -20 /
 \put{\circle{6}} [Bl] at 0  -19.8
 \put{\circle{6}} [Bl] at 60 20.2
 \put{\circle{6}} [Bl] at 60 -19.8
 \thicklines
 \put{\oval(60,20)[t]} [Bl] at 30 20
 \put{\oval(30,10)[b]} [Bl] at 45 -20
 \plot 30 20 0 -20 /
 \endpicture
\end{equation*}

\bigbreak

A \emph{marked diagram} on $2r$ vertices is a graph with $2r$
vertices arranged in two rows of $r$ vertices each, one above the
other, and $r$ edges such that each vertex is incident to
precisely one edge; the rightmost vertex of each `horizontal' edge
(i.e., joining vertices in the same row)  and the bottommost
vertex of each `vertical' edge (i.e., joining vertices in
different rows) may (or may not) be `marked'.

There are $(2r-1)!!$ `unmarked diagrams' and, therefore,
$2^r(2r-1)!!$ marked diagrams. The unmarked diagrams form a basis
of the classical Brauer algebra.

Any such marked diagram represents an element of the centralizer
algebra $\End_{U(V,h)}(\OrV)$. For instance, the marked diagram
$$
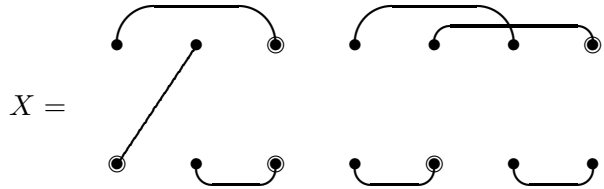
\beginpicture
\setcoordinatesystem units <1.5pt,1.5pt>
 \unitlength=1.5pt 
 \setplotarea x from -50 to 120, y from 30 to 80
 \linethickness=.8pt
 \setplotsymbol ({\sevenrm .})
 \multiput{$\bullet$} at 0 70 20 70 40 70 60 70 /
 \multiput{$\bullet$} at 80 70 100 70 120 70 /
 \multiput{$\bullet$} at 0 40 20 40 40 40 60 40 /
 \multiput{$\bullet$} at 80 40 100 40 120 40 /
 \put{\circle{4}} [Bl] at 40 70.2
 \put{\circle{4}} [Bl] at 0  40.2
 \put{\circle{4}} [Bl] at 120 70.2
 \put{\circle{4}} [Bl] at 40 40.2
 \put{\circle{4}} [Bl] at 80 40.2
 \thicklines
 \put{\oval(40,20)[t]} [Bl] at 20 70
 \put{\oval(40,20)[t]} [Bl] at 80 70
 \put{\oval(40,10)[t]} [Bl] at 100 70
 \put{\oval(20,10)[b]} [Bl] at 30 40
 \put{\oval(20,10)[b]} [Bl] at 70 40
 \put{\oval(20,10)[b]} [Bl] at 110 40
 \plot 20 70 0 40 /
 \put{$X=$} at -20 55
\endpicture
$$
represents the map
\begin{multline*}
\rho(X): v_1\otimes\cdots \otimes v_7\mapsto\\
 \langle v_1\mid Jv_3\rangle\langle v_4\mid v_6\rangle \langle
 v_5\mid Jv_7\rangle
  \sum_{a,b,c=1}^{2n}Jv_2\otimes e_4\otimes Jf_4\otimes e_b\otimes
     Jf_b\otimes e_c\otimes f_c
\end{multline*}

Let $\Dr$ denote the real vector space with a basis formed by the
marked diagrams with $2r$ vertices, numbered from $1$ to $r$ from
left to right in the top row and from $r+1$ to $2r$ from left to
right on the bottom row. The procedure above provides a map from
the set of marked diagrams into $\End_{U(V,h)}(\OrV)$ and hence a
linear map
\begin{equation*}
\rho: \Dr \longrightarrow \End_{U(V,h)}(\OrV).
\end{equation*}
This linear map $\rho$ is onto because of Theorem
\ref{t:centralizador} (or Theorem \ref{t:invariantes}).

\begin{proposition}\label{p:biyeccion}
$\rho$ is a bijection if and only if $n\geq r$.
\end{proposition}

\begin{proof}
Let $X$ be a marked diagram and let us split the edges in $X$
according as wether its rightmost or bottommost vertex is marked
or not:
\begin{equation*}
\begin{split}
X^+&=\{ (p,q)\text{ edge in $X$} \mid p<q\text{ and $q$ is not
marked}\},\\
X^-&=\{ (p,q)\text{ edge in $X$} \mid p<q\text{ and $q$ is
marked}\},
\end{split}
\end{equation*}

Assume $X^+=\{(p_1,q_1),\ldots,(p_s,q_s)\}$ and
$X^-=\{(p_{s+1},q_{s+1}),\ldots,(p_r,q_r)\}$. Let
$\{d_l\}_{l=1}^n$ be an $h$-orthogonal basis of $V$ as a complex
vector space (that is, $h(d_l,d_m)=0$ for $l\ne m$), so that
$\{d_1,\ldots,d_n,Jd_1,\ldots,Jd_n\}$ is an orthogonal basis of
$V$ relative to $\langle\,\mid\,\rangle$. Through the natural
isomorphisms considered in Section 3, $\rho(X)$ corresponds to the
multilinear invariant map:
\begin{equation*}
f_X: (v_1,\ldots,v_{2r})\mapsto
  \pm\hskip -5pt\prod_{(p,q)\in X^+}\hskip -5pt\langle v_p\mid v_q\rangle
   \hskip -2pt\prod_{(p,q)\in X^-}\hskip -5pt\langle v_p\mid Jv_q\rangle
\end{equation*}
(the $\pm$ sign appears due to the skew symmetry of $J$ relative
to $\langle\,\mid\,\rangle$).

If $n\geq r$, take $v_{p_l}=d_l=v_{q_l}$ for $l=1,\ldots,s$ and
$v_{p_l}=Jd_l$, $v_{q_l}=d_l$ for $l=s+1,\ldots,r$. Then
$f_X(v_1,\ldots,v_{2r})\ne 0$, while $f_Y(v_1,\ldots,v_{2r})=0$
for any $Y\ne X$, due to the orthogonality of the chosen basis.
This shows that for $n\geq r$, $\rho$ is one-to-one, and hence a
bijection.

However, if $n\leq r-1$, consider the element $z\in
\End_{U(V,h)}(\OrV)$ given by
\begin{equation*}
z=(1-J_1J_2)(1-J_1J_3)\cdots(1-J_1J_r)\bigl(\sum_{\sigma\in
S_r}(-1)^\sigma\rho(\sigma_r)\bigr),
\end{equation*}
where $(-1)^\sigma$ denotes the signature of $\sigma$. Notice that
\begin{equation}\label{e:z}
z=2^{r-1}e_{\{1,\ldots,r\}}\bigl(\sum_{\sigma\in
S_r}(-1)^\sigma\rho(\sigma_r)\bigr).
\end{equation}

When expanded, $z$ appears as the image under $\rho$ of a
nontrivial linear combination of different marked diagrams without
horizontal edges. For any $\sigma\in S_r$ and any $l\in
\{1,\ldots,r\}$, $\rho(\sigma)J_l=J_{\sigma(l)}\rho(\sigma)$
(remember that $\End_{U(V,h)}(\OrV)$ acts on the right) so, due to
Proposition \ref{p:2.2}(1), for any $\sigma\in S_r$ one has
$e_{\{1,\ldots,r\}}\rho(\sigma)=\rho(\sigma)e_{\{1,\ldots,r\}}$.
Hence the isomorphism
\begin{equation*}
\begin{split}
\Phi_{\{1,\ldots,r\}}:\
\bigl(\OrV\bigr)e_{\{1,\ldots,r\}}\ &\longrightarrow\ \otimes_\C^r V\\
(v_1\otimes\cdots\otimes v_r)e_{\{1,\ldots,r\}}&\mapsto
 v_1\otimes\cdots \otimes v_r
\end{split}
\end{equation*}
given in Proposition \ref{p:2.5}, preserves the action of $S_r$.
Since $n\geq r-1$, $\sum_{\sigma\in S_r}(-1)^\sigma \sigma$ acts
trivially on $\otimes_\C^rV$ and, therefore, $z$ in \eqref{e:z}
acts trivially on $\OrV$. That is, $z=0$. Thus $\rho$ is not
one-to-one in this case.
\end{proof}

\bigbreak

The multiplication (composition of maps) in $\End_{U(V,h)}(\OrV)$
can be lifted to a multiplication in $\Dr$. Let us look at an
example first. Take the following two marked diagrams:
$$
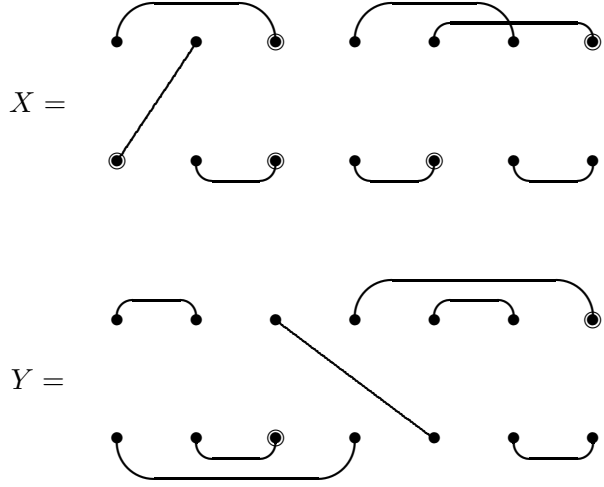
\beginpicture
 \setcoordinatesystem units <1.5pt,1.5pt>
 \unitlength=1.5pt 
 \setplotarea x from -50 to 120, y from -40 to 80
 \linethickness=.8pt
 \setplotsymbol ({\sevenrm .})
 \multiput{$\bullet$} at 0 70 20 70 40 70 60 70 /
 \multiput{$\bullet$} at 80 70 100 70 120 70 /
 \multiput{$\bullet$} at 0 40 20 40 40 40 60 40 /
 \multiput{$\bullet$} at 80 40 100 40 120 40 /
 \multiput{$\bullet$} at 0 0 20 0 40 0 60 0 /
 \multiput{$\bullet$} at 80 0 100 0 120 0 /
 \multiput{$\bullet$} at 0 -30 20 -30 40 -30 60 -30 /
 \multiput{$\bullet$} at 80 -30 100 -30 120 -30 /
 \put{\circle{4}} [Bl] at 40 70.2
 \put{\circle{4}} [Bl] at 0  40.2
 \put{\circle{4}} [Bl] at 120 70.2
 \put{\circle{4}} [Bl] at 40 40.2
 \put{\circle{4}} [Bl] at 80 40.2
 \put{\circle{4}} [Bl] at 120 0.2
 \put{\circle{4}} [Bl] at 40 -29.8
 \thicklines
 \put{\oval(40,20)[t]} [Bl] at 20 70
 \put{\oval(40,20)[t]} [Bl] at 80 70
 \put{\oval(40,10)[t]} [Bl] at 100 70
 \put{\oval(20,10)[b]} [Bl] at 30 40
 \put{\oval(20,10)[b]} [Bl] at 70 40
 \put{\oval(20,10)[b]} [Bl] at 110 40
 \put{\oval(20,10)[t]} [Bl] at 10 0
 \put{\oval(60,20)[t]} [Bl] at 90 0
 \put{\oval(20,10)[t]} [Bl] at 90 0
 \put{\oval(20,10)[b]} [Bl] at 30 -30
 \put{\oval(60,20)[b]} [Bl] at 30 -30
 \put{\oval(20,10)[b]} [Bl] at 110 -30
 \plot 20 70 0 40 /
 \plot 40 0 80 -30 /
 \put{$X=$} at -20 55
 \put{$Y=$} at -20 -15
 \endpicture
$$

Then:
\begin{multline*}
v_1\otimes\cdots\otimes v_7\ \xrightarrow{\ \rho_X\ }\\
 \langle v_1\mid Jv_3\rangle\langle v_4\mid v_6\rangle \langle
  v_5\mid Jv_7\rangle
   \sum_{a,b,c=1}^{2n}Jv_2\otimes e_a\otimes Jf_a\otimes
     e_b\otimes Jf_b\otimes e_c\otimes f_c
\end{multline*}
(with notations already familiar) and
\begin{multline*}
Jv_2\otimes e_a\otimes Jf_a\otimes e_b\otimes Jf_b\otimes
   e_c\otimes f_c\ \xrightarrow{\ \rho_Y\ }\\
    \langle Jv_2\mid e_a\rangle\langle e_b\mid Jf_c\rangle
     \langle Jf_b\mid e_c\rangle
      \sum_{j,k,l=1}^{2n} e_j\otimes e_k\otimes Jf_k\otimes
        f_j\otimes Jf_a\otimes e_l\otimes f_l
\end{multline*}
but $\sum_{a=1}^{2n}\langle Jv_2\mid e_a\rangle Jf_a =
J\bigl(\sum_{a=1}^{2n}\langle Jv_2\mid e_a\rangle f_a\bigr) =
J(Jv_2)=-v_2$, and
\begin{multline*}
\sum_{b,c=1}^{2n}\langle
e_b\mid Jf_c\rangle\langle Jf_b\mid e_c\rangle = -
 \sum_{b,c=1}^{2n}\langle e_b\mid Jf_c\rangle\langle f_b\mid
    Je_c\rangle =\\
   -\sum_{c=1}^{2n}\left\langle \sum_{b=1}^{2n}\langle e_b\mid
    Jf_c\rangle f_b\Bigm| Je_c\right\rangle
 =-\sum_{c=1}^{2n}\langle Jf_a\mid Je_c\rangle
 =-\sum_{c=1}^{2n}\langle e_c\mid f_c\rangle =-2n,
\end{multline*}
 since $J$ is
skew-symmetric relative to $\langle\,\mid\,\rangle$ and $J^2=-1$.
Therefore,
\begin{multline*}
v_1\otimes\cdots\otimes v_7\ \xrightarrow{\,\rho_X\rho_Y\,}\\
  2n \langle v_1\mid Jv_3\rangle\langle v_4\mid v_6\rangle \langle
  v_5\mid Jv_7\rangle
   \sum_{j,k,l=1}^{2n}e_j\otimes e_k\otimes Jf_k\otimes
        f_j\otimes v_2\otimes e_l\otimes f_l
\end{multline*}
which is the image under $\rho$ of $2n$ times the marked diagram
\begin{equation}\label{e:X*Y}
\beginpicture
 \setcoordinatesystem units <1.5pt,1.5pt>
 \unitlength=1.5pt 
 \setplotarea x from -50 to 120, y from -25 to 25
 \linethickness=.8pt
 \setplotsymbol ({\sevenrm .})
 \multiput{$\bullet$} at 0 15 20 15 40 15 60 15 /
 \multiput{$\bullet$} at 80 15 100 15 120 15 /
 \multiput{$\bullet$} at 0 -15 20 -15 40 -15 60 -15 /
 \multiput{$\bullet$} at 80 -15 100 -15 120 -15 /
 \thinlines
 \put{\circle{4}} [Bl] at 40 15.2
 \put{\circle{4}} [Bl] at 120 15.2
 \put{\circle{4}} [Bl] at 40 -14.8
 \thicklines
 \put{\oval(40,20)[t]} [Bl] at 20 15
 \put{\oval(40,20)[t]} [Bl] at 80 15
 \put{\oval(40,10)[t]} [Bl] at 100 15
 \put{\oval(20,10)[b]} [Bl] at 30 -15
 \put{\oval(60,20)[b]} [Bl] at 30 -15
 \put{\oval(20,10)[b]} [Bl] at 110 -15
 \plot 20 15 80 -15 /
 \put{$X*Y=$} at -25 0
 \endpicture
\end{equation}

\bigbreak

Let us consider another example:
$$
\beginpicture
 \setcoordinatesystem units <1pt,1pt>
 \unitlength=1pt 
 \setplotarea x from -30 to 270, y from 0 to 50
 \linethickness=.8pt
 \setplotsymbol ({\sevenrm .})
 \multiput{$\bullet$} at 0 40 20 40 40 40 60 40 80 40 100 40 /
 \multiput{$\bullet$} at 0 10 20 10 40 10 60 10 80 10 100 10 /
 \multiput{$\bullet$} at 160 40 180 40 200 40 220 40 240 40 260 40 /
 \multiput{$\bullet$} at 160 10 180 10 200 10 220 10 240 10 260 10 /
 \put{\circle{6}} [Bl] at 20 10.2
 \put{\circle{6}} [Bl] at 60  10.2
 \put{\circle{6}} [Bl] at 260 40.2
 \thicklines
 \put{\oval(20,10)[t]} [Bl] at 10 40
 \put{\oval(20,10)[t]} [Bl] at 50 40
 \put{\oval(20,10)[t]} [Bl] at 90 40
 \put{\oval(20,10)[b]} [Bl] at 10 10
 \put{\oval(20,10)[b]} [Bl] at 50 10
 \put{\oval(20,10)[b]} [Bl] at 90 10
 \put{\oval(100,20)[t]} [Bl] at 210 40
 \put{\oval(20,10)[t]} [Bl] at 190 40
 \put{\oval(20,10)[t]} [Bl] at 230 40
 \put{\oval(20,10)[b]} [Bl] at 170 10
 \put{\oval(20,10)[b]} [Bl] at 210 10
 \put{\oval(20,10)[b]} [Bl] at 250 10
 \put{$A=$} at -20 25
 \put{,} at 110 25
 \put{$B=$} at 140 25
 \endpicture
$$

Then
\begin{multline*}
v_1\otimes\cdots\otimes v_6\ \xrightarrow{\ \rho_A\ }\\
 \langle v_1\mid v_2\rangle \langle v_3\mid v_4\rangle \langle
   v_5\mid v_6\rangle
   \sum_{a,b,c=1}^{2n}e_a\otimes Jf_a\otimes e_b\otimes
   Jf_b\otimes e_c\otimes f_c
\end{multline*}
while
\begin{multline*}
e_a\otimes Jf_a\otimes e_b\otimes Jf_b\otimes e_c\otimes f_c\
 \xrightarrow{\ \rho_B\ }\\
  \langle e_a\mid Jf_c\rangle \langle Jf_a\mid e_b\rangle
   \langle Jf_b\mid e_c\rangle
    \sum_{j,k,l=1}^{2n} e_j\otimes f_j\otimes e_k\otimes
    f_k\otimes e_l\otimes f_l.
\end{multline*}
Now,
\begin{equation*}
\begin{split}
&\sum_{a,b,c=1}^{2n} \langle e_a \mid Jf_c\rangle
  \langle Jf_a\mid e_b\rangle \langle Jf_b\mid e_c\rangle
  = -\sum_{a,b,c=1}^{2n} \langle e_a \mid Jf_c\rangle
  \langle f_a\mid Je_b\rangle \langle Jf_b\mid e_c\rangle\\
  &\quad = -\sum_{b,c=1}^{2n}\left\langle \sum_{a=1}^{2n}\langle
  e_a\mid Jf_c\rangle f_a\Bigm| Je_b\right\rangle
   \langle Jf_b\mid e_c\rangle
  = -\sum_{b,c=1}^{2n}\langle Jf_c\mid Je_b\rangle
   \langle Jf_b\mid e_c\rangle\\
  &\quad = -\sum_{b,c=1}^{2n}\langle f_c\mid e_b\rangle
   \langle Jf_b\mid e_c\rangle
  =-\sum_{c=1}^{2n}\left\langle J\bigl(\sum_{b=1}^{2n}
   \langle e_b\mid f_c\rangle f_b\bigr)\Bigm| e_c\right\rangle\\
  &\quad
  =-\sum_{c=1}^{2n}\langle Jf_c\mid e_c\rangle =0
\end{split}
\end{equation*}
because of the skew-symmetry of $J$ and since
$\sum_{a=1}^{2n}e_c\otimes f_c=\sum_{a=1}^{2n}f_c\otimes e_c$
(this element of $V\otimes_\R V$ does not depend on the chosen
dual bases). Therefore $\rho_A\rho_B=0$

\medskip

The previous arguments show the general rule to multiply marked
diagrams:

Given two marked diagrams $X$ and $Y$, draw $Y$ below $X$ and
connect the $l^{\text{th}}$ upper vertex of $Y$ with the
$l^{\text{th}}$ lower vertex of $X$, to get a `marked graph'
$G(X,Y)$. For the previously considered marked diagrams $X$ and
$Y$, we have:
$$
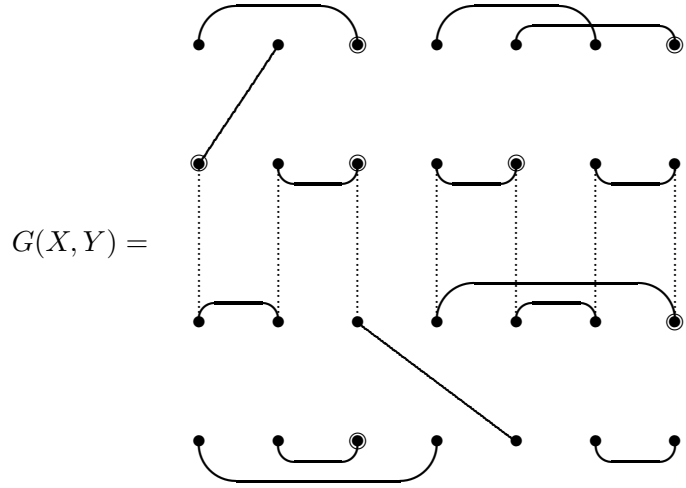
\beginpicture
 \setcoordinatesystem units <1.5pt,1.5pt>
 \unitlength=1.5pt 
 \setplotarea x from -50 to 120, y from -40 to 80
 \linethickness=.8pt
 \setplotsymbol ({\sevenrm .})
 \multiput{$\bullet$} at 0 70 20 70 40 70 60 70 /
 \multiput{$\bullet$} at 80 70 100 70 120 70 /
 \multiput{$\bullet$} at 0 40 20 40 40 40 60 40 /
 \multiput{$\bullet$} at 80 40 100 40 120 40 /
 \multiput{$\bullet$} at 0 0 20 0 40 0 60 0 /
 \multiput{$\bullet$} at 80 0 100 0 120 0 /
 \multiput{$\bullet$} at 0 -30 20 -30 40 -30 60 -30 /
 \multiput{$\bullet$} at 80 -30 100 -30 120 -30 /
 \put{\circle{4}} [Bl] at 40 70.2
 \put{\circle{4}} [Bl] at 0  40.2
 \put{\circle{4}} [Bl] at 120 70.2
 \put{\circle{4}} [Bl] at 40 40.2
 \put{\circle{4}} [Bl] at 80 40.2
 \put{\circle{4}} [Bl] at 120 0.2
 \put{\circle{4}} [Bl] at 40 -29.8
 \thicklines
 \put{\oval(40,20)[t]} [Bl] at 20 70
 \put{\oval(40,20)[t]} [Bl] at 80 70
 \put{\oval(40,10)[t]} [Bl] at 100 70
 \put{\oval(20,10)[b]} [Bl] at 30 40
 \put{\oval(20,10)[b]} [Bl] at 70 40
 \put{\oval(20,10)[b]} [Bl] at 110 40
 \put{\oval(20,10)[t]} [Bl] at 10 0
 \put{\oval(60,20)[t]} [Bl] at 90 0
 \put{\oval(20,10)[t]} [Bl] at 90 0
 \put{\oval(20,10)[b]} [Bl] at 30 -30
 \put{\oval(60,20)[b]} [Bl] at 30 -30
 \put{\oval(20,10)[b]} [Bl] at 110 -30
 \plot 20 70 0 40 /
 \plot 40 0 80 -30 /
 \setdots <2pt>
 \plot 0 40 0 0 /
 \plot 20 40 20 0 /
 \plot 40 40 40 0 /
 \plot 60 40 60 0 /
 \plot 80 40 80 0 /
 \plot 100 40 100 0 /
 \plot 120 40 120 0 /
 \put{$G(X,Y)=$} at -30 20
\endpicture
$$

\bigbreak

 Then
\begin{equation}\label{e:producto}
XY=\gamma(X,Y)\, X*Y
\end{equation}
where $\gamma(X,Y)\in\R$ is defined below and $X*Y$ is the marked
diagram whose vertices are the vertices in the upper row of $X$
and the  vertices in the lower row of $Y$ with the horizontal
edges that appear in these rows. The rightmost vertices of these
horizontal edges inherit the marking in $X$ or $Y$. Moreover,
there is a vertical edge joining any upper vertex of $X$ with a
lower vertex of $Y$ precisely if there is a path in $G(X,Y)$
joining these vertices. The lowermost vertex of any of these
vertical edges is marked if and only if there is an odd number of
marked vertices along the corresponding path in $G(X,Y)$. (For the
example above, $X*Y$ appears in \eqref{e:X*Y}.) Besides,

\begin{enumerate}

\item Let $p$ be any path in $G(X,Y)$, we move the (say) $s$ marks
on the vertices along $p$ to the bottommost vertex and define
\begin{equation*}
\gamma(p)=(-1)^{\text{number of horizontal
moves}}(-1)^{\lfloor\frac{s}{2}\rfloor}
\end{equation*}
($\lfloor x\rfloor$ is the largest integer $\leq x$). Thus, for
instance, consider  the path in the previous example:
$$
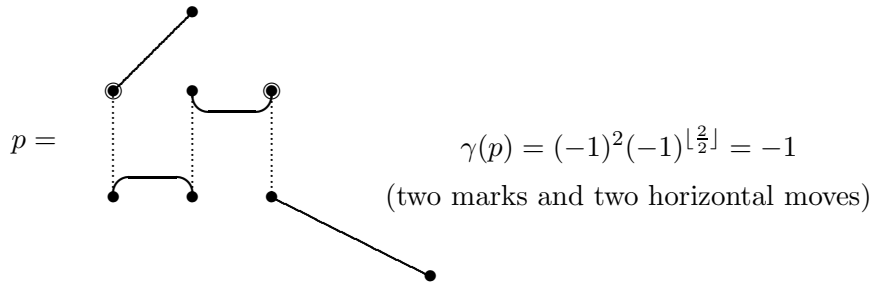
\beginpicture
\setcoordinatesystem units <1.5pt,1pt>
 \unitlength=1.5pt 
 \setplotarea x from -50 to 180, y from -40 to 80
 \linethickness=.8pt
 \setplotsymbol ({\sevenrm .})
 \multiput{$\bullet$} at  20 70 0 40 20 40 40 40 /
 \multiput{$\bullet$} at 0 0 20 0 40 0 80 -30  /
 \put{\circle{4}} [Bl] at 0  40.2
 \put{\circle{4}} [Bl] at 40 40.2
 \thicklines
 \put{\oval(20,10)[b]} [Bl] at 30 40
 \put{\oval(20,10)[t]} [Bl] at 10 0
 \plot 20 70 0 40 /
 \plot 40 0 80 -30 /
 \setdots <2pt>
 \plot 0 40 0 0 /
 \plot 20 40 20 0 /
 \plot 40 40 40 0 /
 \put{$p=$} at -20 20
 \put{$\gamma(p)=(-1)^2(-1)^{\lfloor\frac{2}{2}\rfloor}=-1$} at
    130 20
 \put{(two marks and two horizontal moves)} at 130 0
\endpicture
$$

\item Let $l$ be any loop in $G(X,Y)$, fix any vertex in $l$ and
move all the marks, say $s$, in the vertices of $l$ to this fixed
vertex. Define then
\begin{equation*}
\gamma(l)=\begin{cases}
 0&\text{if $s$ is odd,}\\
 (-1)^{\text{number of horizontal
 moves}}(-1)^{\frac{s}{2}}&\text{if $s$ is even.}
 \end{cases}
\end{equation*}
For instance, taking  the loop of the previous example:
$$
\beginpicture
 \setcoordinatesystem units <1.5pt,1pt>
 \unitlength=1.5pt 
 \setplotarea x from 30 to 230, y from -10 to 50
 \linethickness=.8pt
 \setplotsymbol ({\sevenrm .})
 \multiput{$\bullet$} at 60 40 80 40 100 40 120 40 /
 \multiput{$\bullet$} at 60 0 80 0 100 0 120 0 /
 \put{\circle{4}} [Bl] at 80 40.2
 \put{\circle{4}} [Bl] at 120 0.2
 \thicklines
 \put{\oval(20,10)[b]} [Bl] at 70 40
 \put{\oval(20,10)[b]} [Bl] at 110 40
 \put{\oval(60,20)[t]} [Bl] at 90 0
 \put{\oval(20,10)[t]} [Bl] at 90 0
 \setdots <2pt>
 \plot 60 40 60 0 /
 \plot 80 40 80 0 /
 \plot 100 40 100 0 /
 \plot 120 40 120 0 /
 \put{$l=$} at 40 20
 \put{$\gamma(l)=(-1)^2(-1)^{\frac{2}{2}}=-1$.} at
    200 20
\endpicture
$$

\end{enumerate}

The definitions of $\gamma(p)$ and $\gamma(l)$ are made so as to
take into account the skew-symmetry of $J$ and the fact that
$J^2=-1$.

\medskip

Finally, define
\begin{equation}\label{e:gamma}
\gamma(X,Y)=(2n)^\text{number of loops in $G(X,Y)$}\hskip -5pt
\prod_{\substack{\text{$p$ path}\\ \text{in $G(X,Y)$}}}
 \hskip -10pt\gamma(p)
 \hskip -5pt\prod_{\substack{\text{$l$ loop}\\ \text{in
 $G(X,Y)$}}} \hskip -10pt \gamma(l).
\end{equation}

\bigbreak

The resulting algebra (over the real field) thus defined over
$\Dr$ will be denoted by $\Dr(n)$.

\begin{proposition}\label{p:asociatividad}
$\Dr(n)$ is an associative algebra for any $r,n\in \N$.
\end{proposition}
\begin{proof}
The product in $\Dr(n)$ is defined in such a way as to ensure that
$\rho: \Dr(n)\rightarrow \End_{U(V,h)}(\OrV)$ is a homomorphism of
algebras. Hence the result is obvious for $n\geq r$ by Proposition
\ref{p:biyeccion}. In general, formulas \eqref{e:producto} and
\eqref{e:gamma} show that for $X,Y\in\Dr$,
\begin{equation*}
XY=(2n)^{l(X,Y)}s(X,Y) X*Y,
\end{equation*}
where $l(X,Y)\in \Z_{\geq 0}$ is the number of loops in $G(X,Y)$
and $s(X,Y)\in\{0,1,-1\}$. Both $l(X,Y)$ and $s(X,Y)$ are
independent of $n$. The associativity of $\Dr(n)$ is equivalent to
the validity of $(XY)Z=X(YZ)$ for any marked diagrams $X,Y,Z$, or
to the validity of
\begin{multline*}
(2n)^{l(X,Y)+l(X*Y,Z)}s(X,Y)s(X*Y,Z)\\
=(2n)^{l(Y,Z)+l(X,Y*Z)}s(Y,Z)s(X,Y*Z)
\end{multline*}
which is satisfied if and only if for any marked diagrams $X,Y,Z$
\begin{equation*}
\begin{cases}
 l(X,Y)+l(X*Y,Z)=l(Y,Z)+l(X,Y*Z)&\\
 s(X,Y)s(X*Y,Z)=s(Y,Z)s(X,Y*Z)&
 \end{cases}
\end{equation*}
which does not depend on $n$. Therefore $\Dr(n)$ is associative if
and only if so is $\Dr(m)$ for $m$ large enough, which is indeed
the case.
\end{proof}

\begin{remark}\label{r:drx}
The proof above suggests the consideration of the algebra $\Dr(x)$
over $\R(x)$, with a basis formed by the marked diagrams and
multiplication given by
\begin{equation*}
XY=x^{l(X,Y)}s(X,Y) X*Y
\end{equation*}
for any marked diagrams $X,Y$; in analogy with the Brauer algebras
$\Br_r(x)$ considered in \cite{Wenzl,HW}.
\end{remark}

\begin{remark}\label{r:signed}
Given an edge of a marked diagram, call it \emph{positive} if its
bottommost or rightmost vertex is not marked, and \emph{negative}
otherwise. Hence the marked diagrams can be identified with the
signed diagrams in \cite{Par1,Par2}. The algebra $\Dr(x)$ is then
defined over the same vector space as the Signed Brauer Algebra
defined in these references, although the multiplication is
different.
\end{remark}

\bigbreak

Let us proceed now to give a presentation of $\Dr(x)$ by
generators and relations. We will assume $r\geq 4$, the situation
for $r<3$ is simpler and can be deduced easily along the same
lines. First, let us consider the following marked diagrams

$$
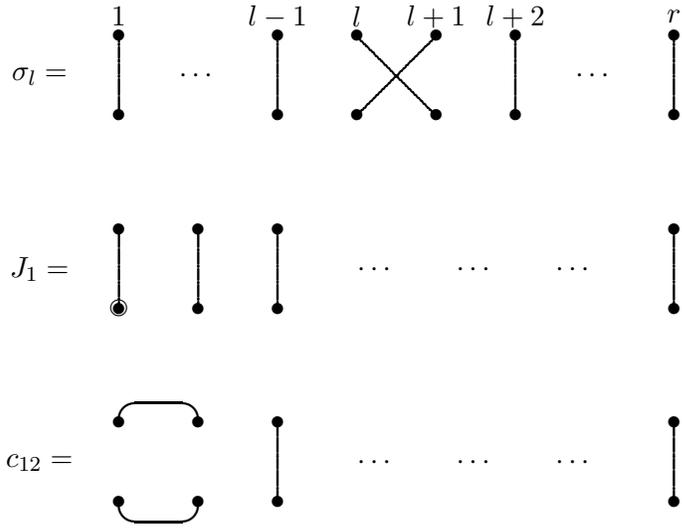
\beginpicture
\setcoordinatesystem units <1.5pt,1.5pt>
 \unitlength=1.5pt 
 \setplotarea x from -30 to 150, y from -20 to 25
 \linethickness=.8pt
 \setplotsymbol ({\sevenrm .})
 \multiput{$\bullet$} at 0 10 40 10 60 10 80 10 100 10 140 10 /
 \multiput{$\bullet$} at 0 -10 40 -10 60 -10 80 -10 100 -10 140 -10 /
 \multiput{$\cdots$} at 20 0 120 0 /
 \put{$1$} at 0 15
 \put{$l-1$} at 40 15
 \put{$l$} at 60 15
 \put{$l+1$} at 80 15
 \put{$l+2$} at 100 15
 \put{$r$} at 140 15
 \plot 0 10 0 -10 /
 \plot 40 10 40 -10 /
 \plot 60 10 80 -10 /
 \plot 80 10 60 -10 /
 \plot 100 10 100 -10 /
 \plot 140 10 140 -10 /
 \put{$\sigma_l=$} at -20 0
\endpicture
$$

$$
\beginpicture
\setcoordinatesystem units <1.5pt,1.5pt>
 \unitlength=1.5pt 
 \setplotarea x from -30 to 150, y from -20 to 25
 \linethickness=.8pt
 \setplotsymbol ({\sevenrm .})
 \multiput{$\bullet$} at 0 10 20 10 40 10  140 10 /
 \multiput{$\bullet$} at 0 -10 20 -10 40 -10  140 -10 /
 \multiput{$\cdots$} at 65 0 90 0 115 0 /
 \plot 0 10 0 -10 /
 \plot 20 10 20 -10 /
 \plot 40 10 40 -10 /
 \plot 140 10 140 -10 /
 \put{\circle{4}} [Bl] at 0 -9.8
 \put{$J_1=$} at -20 0
\endpicture
$$

$$
\beginpicture
\setcoordinatesystem units <1.5pt,1.5pt>
 \unitlength=1.5pt 
 \setplotarea x from -30 to 150, y from -20 to 25
 \linethickness=.8pt
 \setplotsymbol ({\sevenrm .})
 \multiput{$\bullet$} at 0 10 20 10 40 10  140 10 /
 \multiput{$\bullet$} at 0 -10 20 -10 40 -10  140 -10 /
 \multiput{$\cdots$} at 65 0 90 0 115 0 /
 \plot 40 10 40 -10 /
 \plot 140 10 140 -10 /
 \thicklines
 \put{\oval(20,10)[t]} [Bl] at 10 10
 \put{\oval(20,10)[b]} [Bl] at 10 -10
 \put{$c_{12}=$} at -20 0
\endpicture
$$

These generate the algebra $\Dr(x)$, because the $\sigma_l$'s
generate the symmetric group, $\sigma
J_1\sigma^{-1}=J_{\sigma(l)}$ and $\sigma
c_{12}\sigma^{-1}=c_{\sigma(1)\sigma(2)}$ for any $\sigma$ in the
symmetric group.

Then the following relations among these elements are easily
checked:

\begin{enumerate}
\renewcommand{\theenumi}{\roman{enumi}}
\item $\sigma_l^2=1$, $1\leq l\leq r-1$,
\item $\sigma_l\sigma_m=\sigma_m\sigma_l$, $1\leq l,m\leq r-1$,
 $\vert l-m\vert\geq 2$,
\item $(\sigma_l\sigma_{l+1})^3=1$, $1\leq l\leq r-2$,
\item $J_1^2=-1$,
\item $J_1\sigma_l=\sigma_lJ_1$, $2\leq l\leq r-1$,
\item $(J_1\sigma_1)^4=1$,
\item $c_{12}^2=xc_{12}$,
\item $c_{12}\sigma_l=\sigma_lc_{12}$, $3\leq l\leq r-1$,
\item $c_{12}\sigma_1=c_{12}=\sigma_1c_{12}$,
\item $c_{12}\sigma_2c_{12}=c_{12}$,
\item $c_{12}J_1c_{12}=0$,
\item $(\sigma_1J_1+J_1)c_{12}=0=c_{12}(J_1+J_1\sigma_1)$,
\item $\sigma_2\sigma_1J_1\sigma_1\sigma_2c_{12}
  =c_{12}\sigma_2\sigma_1J_1\sigma_1\sigma_2$,
\item $c_{12}\sigma_2\sigma_1\sigma_3\sigma_2
  c_{12}\sigma_2\sigma_3\sigma_1\sigma_2 =
  \sigma_2\sigma_1\sigma_3\sigma_2
  c_{12}\sigma_2\sigma_3\sigma_1\sigma_2c_{12}$.
\end{enumerate}

Notice that (vi) is equivalent to $J_1J_2=J_2J_1$, (vii) to
$J_1c_{12}=-J_2c_{12}$ and $c_{12}J_1=-c_{12}J_2$, (xiii) to
$J_3c_{12}=c_{12}J_3$ and (xiv) to $c_{12}c_{34}=c_{34}c_{12}$.

\medskip

Take the free associative algebra $D$ over $\R(x)$ generated by
elements $\sigma_1,\ldots,\sigma_{r-1},J_1,c_{12}$, subject to the
relations (i)--(xiv) above. $\Dr(x)$ is a quotient of this
algebra, and to show that they are isomorphic it is enough to
check that the dimension of $D$ is $2^r(2r-1)!!$. To do so, first
the subalgebra generated by the $\sigma_l$'s is (isomorphic to)
the group algebra of the symmetric group $S_r$ (in principle it is
a quotient of the group algebra, but the corresponding subalgebra
of $\Dr(x)$ is the whole group algebra). Moreover, define
recursively in $D$ the new elements $J_{l+1}=\sigma_lJ_l\sigma_l$,
$1\leq l\leq r-2$. Because of (v) one has $\sigma
J_l\sigma^{-1}=J_{\sigma(l)}$ for any $\sigma\in S_r\subseteq D$.
Then relation (vi) yields $J_1J_2=J_2J_1$ and with this one proves
easily that $J_lJ_m=J_mJ_l$ for any $l,m$, and that the subalgebra
of $D$ generated by the $\sigma_l$'s and $J_1$ is the span of the
elements $J_\cP\sigma$, where $\cP\subseteq \{1,\ldots,r\}$,
$\sigma\in S_r$ and $J_\cP =\prod_{p\in P}J_p$ ($J_\emptyset =1$).

Now define in $D$ the elements $c_{pq}$ ($p\ne q$) by
$c_{pq}=\sigma c_{12}\sigma^{-1}$, where $\sigma\in S_r$ satisfies
$\sigma(1)=p$, $\sigma(2)=q$. This is well defined by the
relations in (viii) and, because of (ix), $c_{pq}=c_{qp}$ for any
$p,q$. Then relation (xii) is equivalent to $J_3c_{12}=c_{12}J_3$
which, by conjugation with suitable elements of $S_r$, yields
$J_sc_{pq}=c_{pq}J_s$ for different $s,p,q$. Also, relation (xiv)
becomes $c_{12}c_{34}=c_{34}c_{12}$, and again, by conjugation, it
yields $c_{pq}c_{p'q'}=c_{p'q'}c_{pq}$ for different $p,q,p',q'$.
Finally, for distinct elements $p,q,q'$,
$c_{pq}c_{pq'}=c_{pq}(qq')c_{pq}(qq')=c_{pq}(qq')=(qq')c_{pq'}$,
thanks to relation (x) and its conjugates. Besides,
$c_{pq}J_pc_{pq}=0$ by (vi), while
$c_{pq}J_pc_{pq'}=-c_{pq}J_{q'}c_{pq'}=-J_{q'}c_{pq}c_{pq'}$, and
also $c_{pq}J_pc_{pq'}=-c_{pq}J_qc_{pq'}=-c_{pq}c_{pq'}J_q$ for
different $p,q,q'$. With all these relations, any word in the
generators belongs to the linear span of the elements
\begin{equation}\label{e:generadores}
J_{\mathcal Q}c_{p_1q_1}\cdots c_{p_sq_s}J_{\cP^c}\sigma,
\end{equation}
where $p_1<\cdots<p_s,q_1,\ldots,q_s$ are different elements in
$\{1,\ldots,r\}$, ${\mathcal Q}\subseteq \{q_1,\ldots,q_s\}$,
$\cP^c\subseteq \{1,\ldots,r\}\setminus\{p_1,\ldots,p_s\}$ and
$\sigma\in S_r$. But $c_{pq}(pq)=c_{pq}$ by relation (ix) and
$c_{pq}J_q(pq)=c_{pq}(pq)J_p=c_{pq}J_p=-c_{pq}J_q$, so one can
assume that $\sigma(p_l)<\sigma(q_l)$ for any $l=1,\ldots,s$. With
this extra condition on $\sigma$ in \eqref{e:generadores}, each
element in \eqref{e:generadores} is in bijection with a unique
marked diagram.

This finishes the proof of:

\begin{theorem}
Assuming $r\geq 4$, $\Dr(x)$ is the associative algebra over
$\R(x)$ generated by
$\{\sigma_1,\ldots,\sigma_{r-1},J_1,c_{12}\}$, subject to the
relations {\rm (i)--(xiv)}. \qed
\end{theorem}

\begin{remark}
For $r=3$, it is enough to consider relations (i), (iii), (iv),
(vi), (vii), and (ix)--(xiii); while for $r=2$ (i), (iv), (vi),
(vii), (ix), (xi) and (xii) are sufficient.
\end{remark}

\section{Decomposition into irreducibles}

The results in the previous sections, together with
\cite{Benkartetal} (see also \cite{Stembridge}), make it easy to
decompose $\OrV$ into a direct sum of irreducible
$U(V,h)$-modules. First, an element $p\in\{1,\ldots,r\}$ is fixed,
and for simplicity we will take $p=1$. Then, from Proposition
\ref{p:2.2}, Corollary \ref{c:2.3} and Proposition \ref{p:2.5}
\begin{equation*}
\begin{split}
\OrV&=
\hskip -10pt\bigoplus_{1\in\cP\subseteq\{1,\ldots,r\}}
 \hskip -10pt\bigl(\OrV\bigr)e_\cP\\
&\cong \bigoplus_{q=1}^r \binom{r}{q}
 \Bigl(\bigl(\otimes_\C^qV\bigr)
   \otimes_\C\bigl(\otimes_\C^{r-q}V^*\bigr)\Bigr)
\end{split}
\end{equation*}
and, therefore, it is enough to decompose the module (over $\C$)
$V_{q,r-q}:=\bigl(\otimes_\C^qV\bigr)
\otimes_\C\bigl(\otimes_\C^{r-q}V^*\bigr)$ into a direct sum of
irreducible $U(V,h)$-modules. Let us think in terms of the
associated Lie algebra $\uvh$, which is a form of the general
linear Lie algebra $\glv$. The irreducible $\uvh$-submodules of
$V_{q,r-q}$ over $\C$ are exactly the irreducible
$\glv$-submodules, and these are determined in \cite{Stembridge}
and \cite{Benkartetal}: the irreducible $\glv$-submodules of
$V_{q,r-q}$ are in one-to-one correspondence with the pairs
$(\tau,L)$ where:

\begin{enumerate}

\item
 $L=[(m_1,m_1'),\ldots,(m_s,m_s')]$ is a sequence of pairs with
 $1\leq m_1<m_2<\cdots <m_s\leq q$, $m_1',\ldots,m_s'$ are
 different elements in $\{q+1,\ldots,r\}$ ($s\leq\min\{q,r-q\}$).
 $L$ indicates the slots where a contraction is made among $V$ and
 $V^*$.

\item
 $\tau=(\tau^+,\tau^-)$ is a pair of standard rational tableaux,
 where $\tau^+$ (respectively $\tau^-$) is obtained by filling the
 boxes in a Young frame with the numbers in
 $\{1,\ldots,q\}\setminus \{m_1,\ldots,m_s\}$ (resp. in
 $\{q+1,\ldots,r\}\setminus \{m_1',\ldots,m_s'\}$). Being standard
 means that the numbers strictly increase from left to right
 across each row and from top to bottom in each column.

\item
 If $\dim_\C V=n<r$, an extra technical condition has to be
 satisfied (see \cite[Theorem 1.11]{Benkartetal}) that, in
 particular, forces the sum of the number of rows in $\tau^+$ and
 $\tau^-$ to be at most $n$.

\end{enumerate}

\begin{example}
 $r=7$, $q=3$, $L=[(1,6),(2,4)]$, $\tau^+=\begin{Young}
 3\cr\end{Young}$, $\tau^-=\begin{Young} 5\cr 7\cr\end{Young}$.
\end{example}

If we fix a basis of $V$ over $\C$, so that $\glv\cong\gln$, the
complex Lie algebra of $n\times n$ matrices, and consider the
Cartan subalgebra $\h$ formed by the diagonal matrices, let
$\epsilon_l\in\h^*$ ($l=1,\ldots,n$) be given by
$\epsilon_l(\diag(\alpha_1,\ldots,\alpha_n))=\alpha_l$. Then the
highest weight of the irreducible module associated to a pair
$(\tau,L)$ as above is
$\lambda_1\epsilon_1+\cdots+\lambda_n\epsilon_n$ where
$\lambda_1\geq \lambda_2\geq\cdots\geq \lambda_t>0$ are the
lengths of the rows of $\tau^+$, while $-\lambda_n\geq
-\lambda_{n-1}\geq\cdots\geq -\lambda_{n-t'+1}>0$ are the lengths
of the rows of $\tau^-$ and
$\lambda_{t+1}=\cdots=\lambda_{n-t'}=0$.

Once $V_{q,r-q}$ is decomposed into a direct sum of irreducible
modules for $U(V,h)$ over $\C$, what is left to be done is to
check which of these modules remain irreducible as modules for
$U(V,h)$ over $\R$ and which of them do not. The former ones are
the \emph{complex} or \emph{quaternionic} irreducible
representations of $U(V,h)$, while the latter ones are the
\emph{real} representations (notation as in \cite[\S\,26]{FH}. If
$M$ is an irreducible $U(V,h)$-module over $\C$ which is real,
then there exists an irreducible $U(V,h)$-module over $\R$ such
that $M\cong \C\otimes_\R N$. Thus, as a module over $\R$, $M$ is
the direct sum of two copies of $N$.

But if $q\ne r-q$ (in particular, if $r$ is odd), then $i1\in\uvh$
acts as $\bigl(q-(r-q)\bigr)i1=(2q-r)i1\ne 0$ on $V_{q,r-q}=
\bigl(\otimes_\C^qV\bigr)\otimes_\C\bigl(\otimes_\C^{r-q}V^*\bigr)$,
and hence the action of scalar multiplication by imaginary complex
numbers is ``included'' in the action of $\uvh$. Thus all the
irreducible submodules of $V_{q,r-q}$ over $\C$ are complex, so
they are irreducible as modules over $\R$.

The case of $q=r-q$ will be treated in the most interesting case
of $h$ being definite, so $U(V,h)\cong U(n)$ is a compact form of
$GL_\C(V)$. The argument above shows that the action of
$i1\in\uvh$ on $V_{q,q}$ ($r-q=q$) is trivial, so we have to
consider only the action of $\suvh$. Hence, the highest weights of
the irreducible $\glv$-submodules in $V_{q,q}$ are of the form
$\lambda_1\epsilon_1+\cdots+\lambda_n\epsilon_n$, with
$\lambda_1\geq\cdots\geq\lambda_n$ and
$\lambda_1+\cdots+\lambda_n=0$; so that
$\lambda_1\epsilon_1+\cdots+\lambda_n\epsilon_n
=(\lambda_1-\lambda_2)\omega_1+\cdots+
(\lambda_{n-1}-\lambda_n)\omega_{n-1}$, where
$\omega_1=\epsilon_1,\, \omega_2=\epsilon_1+\epsilon_2,\,\ldots\,,
\omega_{n-1}=\epsilon_1+\cdots+\epsilon_{n-1}$, are the
fundamental dominant weights of $\slv$. Notice that the integers
$\lambda_1-\lambda_2,\ldots,\lambda_{n-1}-\lambda_n$, together
with the condition $\lambda_1+\cdots+\lambda_n=0$, determine
$\lambda_1,\ldots,\lambda_n$. The conditions for this highest
weight to yield a real representation are \cite[Proposition
6.24]{FH}:
\begin{align}
&\lambda_l-\lambda_{l+1}=\lambda_{n-l}-\lambda_{n-l+1},\ 1\leq
l\leq n-1,\label{e:realuno}\\
&\text{$n$ odd, or $n=4k$, or $n=4k+2$ and
$\lambda_{2k+1}-\lambda_{2k+2}$ even.}\label{e:realdos}
\end{align}
But \eqref{e:realuno} is equivalent to
$\lambda_l+\lambda_{n+1-l}=0$ for any $l$, which together with the
condition $\lambda_1+\cdots+\lambda_n=0$ yields
$\lambda_l+\lambda_{n+1-l}=0$, which implies
$\lambda_{2k+1}-\lambda_{2k+2}=2\lambda_{2k+1}$ for $n=4k+2$.
Therefore, the condition \eqref{e:realdos} is superfluous, while
the condition \eqref{e:realuno} is equivalent to the restriction
of the Young frames of both $\tau^+$ and $\tau^-$ being the same.

The above discussion is summarized in:

\begin{proposition}\label{p:irreducibles}
Assume $h$ is definite. Then:
\begin{enumerate}

\item For $p\ne q$, the irreducible $U(V,h)$-submodules of
$\bigl(\otimes_\C^pV\bigr)\otimes_\C\bigl(\otimes_\C^qV^*\bigr)$
over $\C$ are all complex, so they remain irreducible as modules
over $\R$.

\item The same happens if $p=q$ for the irreducible
$U(V,h)$-modules of
$\bigl(\otimes_\C^qV\bigr)\otimes_\C\bigl(\otimes_\C^qV^*\bigr)$
over $\C$ which correspond to pairs $(\tau,L)$ where the Young
frames of $\tau^+$ and $\tau^-$ are different.

\item The irreducible $U(V,h)$-submodules of
$\bigl(\otimes_\C^pV\bigr)\otimes_\C\bigl(\otimes_\C^qV^*\bigr)$
over $\C$ which correspond to pairs $(\tau,L)$, with equal Young
frames of $\tau^+$ and $\tau^-$ are real, so they split into a
direct sum of two copies of an irreducible $U(V,h)$-module over
$\R$.

\end{enumerate}
\end{proposition}

As a first example, consider $V\otimes_\R V$, which splits as:
\begin{equation*}
\begin{split}
V\otimes_\R V&=(V\otimes_\R V)e_1\oplus (V\otimes_\R V)e_2\simeq
    (V\otimes_\C V)\oplus (V\otimes_\C V^*)\\
    &\simeq S^2(V)\oplus \Lambda^2(V)\oplus \C\oplus \slv
\end{split}
\end{equation*}
where $e_1=\frac{1}{2}(1-J_1J_2)$, $e_2=\frac{1}{2}(1+J_1J_2)$.
Here $S^2(V)$ (the symmetric tensors) and $\Lambda^2(V)$
(skew-symmetric tensors) are irreducible as $U(V,h)$-modules over
$\R$, while $\C$ is a direct sum of two trivial one-dimensional
modules over $\R$, and $\slv=\suvh\oplus i\suvh$ is a direct sum
of two copies of $\suvh$.

\bigbreak

In the remaining part of this paper, we will consider the
motivating example of Gray and Hervella \cite{Gray-Hervella}
considered in the Introduction, as well as another related example
by Abbena and Garbiero \cite{Abbena-Garbiero}.

\begin{example} {\bf (Gray-Hervella 1978 \cite{Gray-Hervella})}

Since $V\cong V^*$ as modules for $U(V,h)$ over $\R$, the problem
described in the Introduction amounts to decompose
\begin{equation*}
\begin{split}
W&=\{ x\in \otimes_\R^3 V : x(23)=-x=xJ_2J_3\}\\
&=\{x\in (\otimes_\R^3V)\frac{1}{2}(1-J_2J_3) : x(23)=-x\}
\end{split}
\end{equation*}
into a direct sum of irreducible submodules for $U(V,h)$. First
notice that
\begin{equation*}
\begin{split}
\frac{1}{2}(1-J_2J_3)
  &=\frac{1}{2}(1-J_2J_3)\frac{1}{2}(1-J_1J_3+1+J_1J_3)\\
  &=e_{\{1,2,3\}} +e_{\{2,3\}}
\end{split}
\end{equation*}
(a sum of orthogonal idempotents) in the notation of Section 2.
Hence, by Proposition \ref{p:2.5}
\begin{equation*}
\begin{split}
\bigl(\otimes_\R^3V\bigr)\frac{1}{2}(1-J_2J_3)
  &=\bigl(\otimes_\R^3V\bigr)e_{\{1,2,3\}} \oplus
     \bigl(\otimes_\R^3V\bigr)e_{\{2,3\}}\\
  &\cong (V\otimes_\C V\otimes_\C V)\oplus
      (V^*\otimes_\C V\otimes_\C V)
\end{split}
\end{equation*}
and from this isomorphism, it immediately follows that
\begin{equation*}
W\cong \bigl(V\otimes_\C\Lambda^2(V)\bigr)\oplus
  \bigl(V^*\otimes_\C\Lambda^2(V)\bigr),
\end{equation*}
and it is enough to decompose each one of these two summands into
irreducible $\glv$-modules.

The first summand is $(V\otimes_\C V\otimes_\C
V)\frac{1}{2}(1-(23))$ and, since in the group algebra $\C S_3$
the idempotent $\frac{1}{2}(1-(23))=e_{T_1}+e_{T_2}$ is the sum of
two orthogonal primitive idempotents, where
\begin{equation*}
T_1=\begin{Young} 2\cr 3\cr 1\cr\end{Young},\qquad
T_2=\begin{Young} 2&1\cr 3\cr\end{Young},
\end{equation*}
(that is, $e_{T_1}=\frac{1}{6}\sum_{\sigma\in
S_3}(-1)^\sigma\sigma$ and
$e_{T_2}=\frac{1}{3}\bigl(1+(12)\bigr)\bigl(1-(23)\bigr)$), it
follows that $V\otimes_\C\Lambda^2(V) =
(\otimes_\C^3V)e_{T_1}\oplus (\otimes_\C^3V)e_{T_2}$ (direct sum
of two irreducible modules if $\dim_\C V\geq 3$), which under the
isomorphisms correspond to:
\begin{equation*}
\begin{split}
W_1&=\{x\in \otimes_\R^3V : xJ_1=xJ_2=xJ_3,\, xe_{T_1}=x\},\\
W_2&=\{x\in \otimes_\R^3V : xJ_1=xJ_2=xJ_3,\, xe_{T_2}=x\}.
\end{split}
\end{equation*}
But $xe_{T_2}=x$ if and only if $x(23)=-x$ and
$x(1+(123)+(132))=0$. For any $x$ satisfing these two conditions,
one checks easily that the condition $xJ_1=xJ_2$ follows from
$xJ_2=xJ_3$. Hence
\begin{equation*}
W_2=\{ x\in W : x(1+(123)+(132))=0\},
\end{equation*}
and, similarly,
\begin{equation*}
W_1=\{ x\in W : x(12)=-x\}.
\end{equation*}
On the other hand, assuming $\dim_\C V\geq 3$,
$V^*\otimes_\C\Lambda^2(V)$ decomposes \cite{Benkartetal} into:

\begin{itemize}

\item $\ker c$, where $c:V^*\otimes_\C\Lambda^2(V)\rightarrow V$,
$f\otimes(u\otimes v-v\otimes u)\mapsto f(u)v-f(v)u$, which
corresponds to the pair $(\tau,L)$, with $\tau=\Bigl(\,
 \begin{Young} 1\cr 2\cr\end{Young},\begin{Young}
 3\cr\end{Young}\Bigr)$
 and $L=\emptyset$.

\item $\{\sum_{l=1}^n u_l^*\otimes(u_l\otimes v-v\otimes u_l) :
v\in V\}\cong V$, where $\{u_l\}_{l=1}^n$ and $\{u_l^*\}_{l=1}^n$
are dual bases of $V$ and $V^*$ over $\C$. This is a `diagonal'
submodule of the ones that correspond to the pairs $(\tau,L)$ with
$L=[(1,3)]$ and $\tau=(\boxed{2},\emptyset)$ and with $L=[(2,3)]$
and $\tau=(\boxed{1},\emptyset)$.

\end{itemize}

Under the isomorphisms, these submodules correspond to:
\begin{equation*}
\begin{split}
W_3&=\{ x\in W : xJ_1=-xJ_2\text{ and }xc_{12}=0\}\\
W_4&=\{ \sum_{l=1}^{2n}(e_l\otimes f_l\otimes
v)\frac{1}{4}(1+J_1J_3)(1-(23) : v\in V\}
\end{split}
\end{equation*}
where the $e_l$'s and $f_l$'s constitute dual bases of $V$ over
$\R$ relative to $\langle\,\mid\,\rangle$.

The situation for $\dim_\C V\leq 2$ is simpler.

We recover in this way the decomposition given in
\cite{Gray-Hervella}.

\end{example}

\begin{example} {\bf (Abbena-Garbiero 1988 \cite{Abbena-Garbiero})}

One has to decompose
\begin{equation*}
\begin{split}
K&=\{ x\in \otimes_\R^3V : x(23)=-x=-xJ_2J_3\}\\
 &=\{ x\in (\otimes_\R^3V)\frac{1}{2}(1+J_2J_3) : x(23)=-x\} .
\end{split}
\end{equation*}
As before, $\frac{1}{2}(1+J_2J_3)=e_{\{1,2\}} +e_{\{1,3\}}$, so
\begin{equation*}
(\otimes_\R^3V)\frac{1}{2}(1+J_2J_3)\cong
 (V\otimes_\C V\otimes_\C V^*) \oplus
   (V\otimes_\C V^*\otimes_\C V)
\end{equation*}
by means of the isomorphism $\Phi$ given by
\begin{equation*}
\Phi\bigl( (v_1\otimes v_2\otimes v_3)\frac{1}{2}(1+J_2J_3)\bigr)
 = v_1\otimes v_2\otimes h(-,v_3) + v_1\otimes h(-,v_2)\otimes
 v_3.
\end{equation*}
Now, the following diagram is commutative:
\begin{equation*}
\begin{CD}
\bigl(\otimes_\R^3V\bigr)\frac{1}{2}(1+J_2J_3) @>\Phi>>
  (V\otimes_\C V\otimes_\C V^*) \oplus
   (V\otimes_\C V^*\otimes_\C V)\\
  @V(23)VV  @VV{\text{flip}}V \\
\bigl(\otimes_\R^3V\bigr)\frac{1}{2}(1+J_2J_3) @>\Phi>>
  (V\otimes_\C V\otimes_\C V^*) \oplus
   (V\otimes_\C V^*\otimes_\C V)
\end{CD}
\end{equation*}
Here we have on the left the right action of the transposition
$(23)$, while on the right $(v_1\otimes v_2\otimes
f)\text{flip}=v_1\otimes f\otimes v_2$ and $(v_1\otimes f\otimes
v_2)\text{flip}=v_1\otimes v_2\otimes f$ for any $v_1,v_2\in v$
and $f\in V^*$.

Therefore, the linear map given by
\begin{equation*}
\begin{split}
K=\bigl(\otimes_\R^3 V\bigr)\frac{1}{4}(1+J_2J_3)(1-(23))\ &
  \longrightarrow\  V\otimes_\C\otimes_\C V^*\\
(v_1\otimes v_2\otimes v_3)\frac{1}{4}(1+J_2J_3)(1-(23)) &
 \mapsto \frac{1}{2}\bigl(v_1\otimes v_2\otimes h(-,v_3) -\\
 &\hskip 1in
    v_1\otimes v_3-h(-,v_2)\bigr)
\end{split}
\end{equation*}
is an isomorphism of $U(V,h)$-modules.

If $\dim_\C V\geq 3$, $V\otimes_\C V\otimes_\C V^*$ decomposes
into the direct sum of the irreducible $\glv$-modules
corresponding to the pairs $(\tau,L)$ in the list (see
\cite{Benkartetal}):
\begin{equation*}
\begin{split}
&\tau=\bigl(\, \begin{Young} 1&2\cr\end{Young},
  \begin{Young} 3\cr\end{Young}\bigr),\quad L=\emptyset,\\
&\tau=\Bigl(\, \begin{Young} 1\cr 2\cr\end{Young},
  \begin{Young} 3\cr\end{Young}\Bigr),\quad L=\emptyset,\\
&\tau=\bigl(\, \begin{Young} 1\cr\end{Young},
  \emptyset\bigr),\quad L=[(2,3)],\\
&\tau=\bigl( \begin{Young} 2\cr\end{Young},
  \emptyset\bigr),\quad L=[(1,3)].
\end{split}
\end{equation*}

With $c: V\otimes_\C V\otimes_\C V^*\rightarrow V$, $v_1\otimes
v_2\otimes f\mapsto f(v_1)v_2$, the first two modules correspond
to $\bigl(S^2(V)\otimes_\C V^*\bigr)\cap \ker c$ and
$\bigl(\Lambda^2(V)\otimes_\C V^*\bigr)\cap \ker c$, while the
last two modules are isomorphic to $V$. One recovers from here the
decomposition given in \cite{Abbena-Garbiero}. The details are
left to the reader.

\end{example}

\begin{acknow} The author is indebted to  Pedro
Mart{\'\i}nez Gadea, who introduced him to the problems in
\cite{Gray-Hervella,Abbena-Garbiero} and sent a copy of
\cite{Pedroetal}, and to Georgia Benkart for very illuminating
conversations about this work.
\end{acknow}


\providecommand{\bysame}{\leavevmode\hbox
to3em{\hrulefill}\thinspace}
\providecommand{\MR}{\relax\ifhmode\unskip\space\fi MR }

\end{document}